\documentclass[11pt,twocolumn]{article}
\usepackage{amssymb,amsmath, amsthm}
\usepackage{amssymb}
\usepackage{framed}
\usepackage{mathrsfs}
\usepackage{abstract}

\usepackage{enumerate}
\usepackage{tikz}
\usetikzlibrary{matrix}
\usepackage[all]{xy}
\usetikzlibrary{shapes}
\usetikzlibrary{decorations.markings}
\usepackage{hyperref}
\usepackage{multirow}
\usepackage{bm}
\usepackage[algoruled,boxruled, shortend,noline]{algorithm2e}
\usepackage{algpseudocode}
\usepackage{cite}
\usepackage{enumitem}
\setenumerate[0]{leftmargin=*}

\numberwithin{equation}{section}

\newcommand{\inclu}[0] {\ar@{^{(}->}}

\newcommand{\dist}{{\rm dist}}
\newcommand{\prox}{{\rm prox}}
\newcommand{\proj}{{\rm proj}}
\newcommand{\R}{{\mathbb R}}

\newlength{\bibitemsep}
\newlength{\bibparskip}
\let\oldthebibliography\thebibliography
\renewcommand\thebibliography[1]{%
  \oldthebibliography{#1}%
  \setlength{\parskip}{\bibitemsep}%
  \setlength{\itemsep}{\bibparskip}%
}

\usepackage[margin=1in]{geometry}

\newcommand{\argmin}{\operatornamewithlimits{argmin}}

\SetAlgoSkip{bigskip}
\SetKwInput{Input}{Input\hspace{34pt}{ }}
\SetKwInput{Output}{Output\hspace{27pt}{ }}
\SetKwInput{Initialize}{Initialize\hspace{2pt}{ }}
\SetKwInput{Parameters}{Parameters\hspace{2pt}{ }}

\newtheorem{theorem}{Theorem}[section]

\newtheorem{defn}[theorem]{Definition}

\newtheorem{example}{Example}[section]

\title{The proximal point method revisited\footnote{Submitted to SIAG/OPT Views and News.}}

\begin{document}
	
	\author{Dmitriy Drusvyatskiy 
		\thanks{University of Washington, Department of Mathematics, 
			Seattle, WA 98195;
				\texttt{www.math.washington.edu/{\raise.17ex\hbox{$\scriptstyle\sim$}}ddrusv/}}}	

\date{}

%


\twocolumn[
\maketitle 
\begin{onecolabstract}
	In this short survey, I revisit the role of the proximal point method in large scale optimization. I focus on three recent examples: a proximally guided subgradient method for weakly convex stochastic approximation, the prox-linear algorithm for minimizing compositions of convex functions and smooth maps, and Catalyst generic acceleration  for regularized Empirical Risk Minimization.
\end{onecolabstract}
\bigskip
]

\section{Introduction}
The proximal point method is a conceptually simple   algorithm for minimizing a function $f$ on $\R^d$.  Given an iterate $x_t$, the method  defines  $x_{t+1}$ to be any minimizer of the proximal subproblem
$$\argmin_{x}~\left\{f(x)+\tfrac{1}{2\nu}\|x-x_t\|^2\right\},$$
for an appropriately chosen parameter $\nu>0$. At first glance, each proximal subproblem seems no easier than minimizing $f$ in the first place. On the contrary, the addition of the quadratic penalty term often regularizes the proximal subproblems and makes them well conditioned. Case in point, the  subproblem may become convex despite $f$ not being convex; and even if $f$ were convex, the subproblem has a larger strong convexity parameter thereby facilitating faster numerical methods. 

Despite the improved conditioning, each proximal subproblem still requires invoking an iterative solver. For this reason, the proximal point method has predominantly been thought of as a theoretical/conceptual algorithm, only guiding algorithm design and analysis rather than being implemented directly. One good example is the proximal bundle method \cite{bundle}, which approximates each proximal subproblem by a cutting plane model. In the past few years, this viewpoint has undergone a major revision. In a variety of circumstances, the proximal point method (or a close variant) with a judicious choice of the control parameter $\nu>0$ and an appropriate iterative method for the subproblems can lead to practical and theoretically sound numerical methods. 
In this article, I will briefly describe three recent examples of this trend: 
\begin{itemize}
	\item a subgradient method for weakly convex stochastic approximation problems \cite{prox_guide_subgrad},
	\item  the prox-linear algorithm for minimizing compositions of convex functions and smooth maps \cite{prox_error,prox_lin_paq,quad_conv,nest_GN,prox,composite_cart},
	\item Catalyst generic acceleration schema \cite{catalyst} for regularized Empirical Risk Minimization.
\end{itemize}

In this article, I will focus only on the proximal point method for minimizing functions, as outlined above. The proximal point methodology applies much more broadly to monotone operator inclusions; I refer the reader to the monograph of Bauschke and Combette \cite{BC_book} or the seminal work of Rockafellar \cite{mon_rock}.

\section{Notation}
The following two constructions will play a basic role in the article. For any closed function $f$ on $\R^d$,   the {\em Moreau envelope} and the {\em proximal map} are 
\begin{align*}
f_{\nu}(z)&:=\inf_{x}~\left\{f(x)+\tfrac{1}{2\nu}\|x-z\|^2\right\},\\
\prox_{\nu f}(z)&:=\argmin_{x}~\left\{f(x)+\tfrac{1}{2\nu}\|x-z\|^2\right\},
\end{align*}
respectively.
In this notation, the proximal point method is simply the fixed-point recurrence on the proximal map:\footnote{To ensure that $\prox_{\nu f}(\cdot)$ is nonempty, it suffices to assume that  $f$ is bounded from below.}  $${\bf Step\, }t: \qquad \textrm{choose }x_{t+1}\in \prox_{\nu f}(x_t).$$

Clearly, in order to have any hope of solving the proximal subproblems, one must ensure that they are convex. Consequently, the class of weakly convex functions forms the natural setting for the proximal point method. 
\begin{defn}{\rm
A function $f$ is called {\em $\rho$-weakly convex} if the assignment $x\mapsto f(x)+\frac{\rho}{2}\|x\|^2$ is a convex function.}
\end{defn}

 For example, a $C^1$-smooth function with $\rho$-Lipschitz gradient is $\rho$-weakly convex, while a $C^2$-smooth function $f$ is 
$\rho$-weakly convex precisely when the minimal eigenvalue of its Hessian is uniformly bounded below by $-\rho$.
In essence, weak convexity precludes functions that have downward kinks. For instance, $f(x):=-\|x\|$ is not weakly convex since no addition of a quadratic makes the resulting function convex.

 Whenever $f$ is $\rho$-weakly convex and the proximal parameter  $\nu$ satisfies  $\nu<\rho^{-1}$, each proximal subproblem is itself convex and therefore globally tractable. Moreover, in this setting, the Moreau envelope is $C^1$-smooth with the gradient
\begin{equation}\label{eqn:grad_form}
\nabla f_{\nu}(x)=\nu^{-1}(x-\prox_{\nu f}(x)).
\end{equation}
 Rearranging the gradient formula yields the useful interpretation of the proximal point method as gradient descent on the Moreau envelope
$$x_{t+1}=x_t-\nu\nabla f_{\nu}(x_t).$$

In summary, the Moreau envelope $f_{\nu}$ serves as a $C^1$-smooth approximation of $f$ for all small $\nu$. Moreover, the two conditions $$\|\nabla f_{\nu}(x_{t})\|< \varepsilon$$ and
$$\|\nu^{-1}(x_t-x_{t+1})\|<\varepsilon,$$
are equivalent for the proximal point sequence $\{x_t\}$.
  Hence, the step-size $\|x_t-x_{t+1}\|$ of the proximal point method serves as a convenient termination criteria.

\subsection{Examples of weakly convex functions}
Weakly convex functions are widespread in applications and are typically easy to recognize. One common source of weakly convex functions is the composite problem class:
\begin{equation}\label{eqn:comp}
\min_{x}~ F(x):=g(x)+h(c(x)),
\end{equation}
where $g\colon \R^d\to\R\cup\{+\infty\}$ is a closed convex function, $h\colon\R^m\to\R$ is convex and $L$-Lipschitz, and $c\colon\R^d\to\R^m$ is a $C^1$-smooth map with $\beta$-Lipschitz gradient. An easy argument shows that  $F$ is $L\beta$-weakly convex. This is a worst case estimate. In concrete circumstances, the composite function $F$ may have a much more favorable weak convexity constant (e.g., phase retrieval  \cite[Section 3.2]{duchi_ruan_PR}).


\begin{example}[Additive composite]\label{exa:add_comp}
	{\rm The most prevalent example is additive composite minimization. In this case, the map $c$ maps to the real line and $h$ is the identity function:
		\begin{equation}\label{eqn:add_comp}
		\min_{x}~ c(x)+g(x).
		\end{equation}
		Such problems appear often in statistical learning and imaging. A variety of specialized algorithms are available; see for example Beck and Teboulle \cite{smoothing_beckT} or Nesterov \cite{nest_conv_comp}.

	}
\end{example}

\begin{example}[Nonlinear least squares]\label{exa:nls}
	{\rm
		The composite problem class also captures nonlinear least squares problems with bound constraints:
		\begin{align*}
		\min_x~ \|c(x)\|_2\quad \textrm{subject to}\quad l_i\leq x_i\leq u_i ~\forall i.
		\end{align*}
		Such problems pervade engineering and scientific applications.

	}
\end{example}

\begin{example}[Exact penalty formulations]\label{exa:ep}
	{\rm
		Consider a nonlinear optimization problem:
		\begin{align*}
		\min_x~ \{f(x): G(x)\in \mathcal{K}\},
		\end{align*}
		where $f$ and $G$ are smooth maps and 
		$\mathcal{K}$ is a closed convex cone.
		An accompanying {\em penalty formulation} -- ubiquitous in nonlinear optimization 
		-- takes the form 
		$$\min_x~ f(x)+\lambda \cdot \dist_{\mathcal{K}}(G(x)),$$
		where $\dist_{\mathcal{K}}(\cdot)$ is the  distance  to $\mathcal{K}$ in some norm. 
		Historically, exact penalty formulations served as the  early motivation for the  class \eqref{eqn:comp}.	
	}
\end{example}

\begin{example}[Robust phase retrieval]
	{\rm
	Phase retrieval is a common computational problem, with applications in diverse areas, such as imaging, X-ray crystallography, and speech processing. For simplicity, I will focus on the version of the problem over the reals.
The (real) phase retrieval problem seeks to determine a point $x$ satisfying the magnitude conditions, $$|\langle a_i,x\rangle|\approx b_i\quad \textrm{for }i=1,\ldots,m,$$ where $a_i\in \R^d$ and $b_i\in\R$ are given. Whenever there are gross outliers in the measurements $b_i$, the following robust formulation of the problem is appealing \cite{eM,duchi_ruan_PR, proj_weak_dim}:
$$\min_x ~\tfrac{1}{m}\sum_{i=1}^m |\langle a_i,x\rangle^2-b_i^2|.$$
Clearly, this is an instance of \eqref{eqn:comp}. For some recent perspectives on phase retrieval, see the survey \cite{luke_news_views}. There are numerous recent  nonconvex approaches to phase retrieval, which rely on alternate problem formulations; e.g., \cite{wirt_flow,rand_quad,phase_nonconv}.}	
\end{example}

\begin{example}[Robust PCA]
{\rm
In robust principal component analysis, one seeks to identify sparse corruptions of a low-rank matrix \cite{rob_cand,chand}. One typical example is image deconvolution,  where the low-rank structure  models the background of an image while the sparse corruption models the foreground. 
 Formally, given a $m\times n$ matrix $M$, the goal is to find a decomposition $M=L+S$, where $L$ is low-rank and $S$ is sparse. A common formulation of the problem reads:
$$\min_{U\in \R^{m\times r},V\in \R^{n\times r}}~ \|UV^T-M\|_1,$$ 
where $r$ is the target rank. }
\end{example}


\begin{example}[Censored $\mathbb{Z}_2$ synchronization]
	{\rm
A synchronization problem over a graph is to estimate group elements $g_1,\ldots, g_n$ from  pairwise products  $g_ig_j^{-1}$ over a set of edges $ij\in E$. For a list of application of such problem see \cite{ban_boum,ang_sing,abbe_band}, and references therein. A simple instance is $\mathbb{Z}_2$ synchronization, corresponding to the group on two elements $\{-1,+1\}$. The popular problem of detecting communities in a network, within the Binary Stochastic Block Model (SBM), can be modeled using $\mathbb{Z}_2$ synchronization. 

Formally, given a partially observed matrix $M$, the goal is to recover a vector $ \theta\in \{\pm 1\}^d$, satisfying $M_{ij}\approx \theta_i \theta_j$ for all $ij\in E$. When the entries of $M$ are corrupted by adversarial sign flips, one can postulate the following formulation
$$\min_{\theta\in \R^{d}}~ \|P_{E}(\theta\theta^T-M)\|_1,$$ 
where the operator $P_E$ records the entries indexed by the edge set $E$. Clearly, this is again an instance of \eqref{eqn:comp}.
}
\end{example}

\section{The proximally guided subgradient method}
As the first example of contemporary applications of the 
proximal point method,  consider the problem of minimizing the expectation:\footnote{For simplicity of the exposition, the minimization problem is unconstrained. Simple constraints can be accommodated using a projection operation.}
$$\min_{x\in \R^d}~ F(x)=\mathbb{E}_{\zeta} f(x,\zeta).$$
Here, $\zeta$ is a random variable, and the only access  to $F$ is by sampling $\zeta$.
 It is difficult to overstate the importance of this problem class (often called {\em stochastic approximation}) in large-scale optimization; see e.g. \cite{BB,jordan}.

When the problem is convex, the stochastic subgradient method \cite{stochave,rob_mon,latest_subgrad} has strong theoretical guarantees and is often the method of choice. 
In contrast, when applied to nonsmooth and nonconvex problems, the behavior of the method is poorly understood. The recent paper  \cite{prox_guide_subgrad} shows how to use the proximal point method to guide the subgradient iterates in this broader setting, with rigorous guarantees.

Henceforth, assume that the function $x\mapsto f(x,\zeta)$ is $\rho$-weakly convex and $L$-Lipschitz for each $\zeta$. Davis and Grimmer \cite{prox_guide_subgrad} proposed the scheme outlined in Algorithm~\ref{alg:proxguide}.

\begin{algorithm}
	\KwData{$x_0\in \R^d$, $\{j_t\}\subset\mathbb{N}$, $\{\alpha_j\}\subset\R_{++}$}
	\For{t=0,\ldots,T}{
	Set $y_0=x_t$\;
	\For{$j=0,\ldots,j_t-2$}{
		Sample $\zeta$ and choose $v_j\in\partial (f(\cdot,\zeta)+\rho\|\cdot-x_t\|^2)(y_j)$\;
		$y_{j+1}= y_j-\alpha_jv_j$
}
$x_{t+1}= \frac{1}{j_t}\sum_{j=0}^{j_t-1}y_j$}
\caption{Proximally guided stochastic subgradient method}\label{alg:proxguide}
\end{algorithm}

The method proceeds by applying a proximal point method with each subproblem approximately solved by a stochastic subgradient method. The intuition is that each proximal subproblem is $\rho/2$-strongly convex and therefore according to well-known results (e.g. \cite{a_simp_app,Rakhlin_subgrad,hazan_subgrad,MR3353214}), the stochastic subgradient method should converge at the rate $O(\frac{1}{T})$  on the subproblem, in expectation. This intuition is not quite correct because the objective function of the subproblem is not globally Lipschitz -- a key assumption for the $O(\frac{1}{T})$ rate. Nonetheless, the authors show that warm-starting the subgradient method for each proximal subproblem with the current proximal iterate corrects this issue, yielding a favorable guarantees \cite[Theorem 1]{prox_guide_subgrad}. 

To describe the rate of convergence, set
$j_t=t+\lceil 648\log(648)\rceil$ and $\alpha_j=\tfrac{2}{\rho(j+49)}$ in Algorithm~\ref{alg:proxguide}. Then the scheme will generate an iterate $x$ satisfying 
$$\mathbb{E}_{\zeta}[\|\nabla F_{2\rho}(x)\|^2]\leq \varepsilon$$
after at most 
$$O\left(\frac{\rho^2(F(x_0)-\inf  F)^2}{\varepsilon^2}+\frac{L^4 \log^{4}(\varepsilon^{-1})}{\varepsilon^2}\right)$$
subgradient evaluations. This rate agrees with analogous guarantees for stochastic gradient methods for smooth nonconvex functions \cite{gl_stoch}. 
It is also worth noting that convex constraints on $x$ can be easily incorporated into Algorithm~\ref{alg:proxguide} by introducing a nearest-point projection in the definition of $y_{j+1}$. 

\section{The prox-linear algorithm}	
For well-structured weakly convex problems, one can hope for faster numerical methods than the subgradient scheme. In this section, I will focus on the composite problem class \eqref{eqn:comp}. To simplify the exposition, I will assume $L=1$, which can always be arranged by rescaling.

Since composite functions are weakly convex, one could apply the proximal point method directly, while setting the parameter $\nu\leq\beta^{-1}$. Even though the proximal subproblems are strongly convex, they are not in a form that is most amenable to convex optimization techniques. Indeed, most convex optimization algorithms are designed for minimizing a sum of a convex function and a composition of a convex function with a {\em linear} map. This observation suggests introducing the following modification to the proximal-point algorithm. Given a current iterate $x_t$, the {\em prox-linear method} sets
\begin{align*}
x_{t+1}=\argmin_x \{F(x;x_t)+\tfrac{\beta}{2}\|x-x_t\|^2\},
\end{align*}
where $F(x;y)$ is the local convex model 
$$F(x;y):=g(x)+h\left(c(y)+\nabla c(y)(x-y)\right).$$
In other words, each proximal subproblem is approximated by linearizing the smooth map $c$ at the current iterate $x_t$.

The main advantage is that each subproblem is now a sum of a strongly convex function and a composition of a Lipschitz convex function with a linear map. A variety of methods utilizing this structure can be formally applied; e.g. smoothing \cite{smooth_min_nonsmooth}, saddle-point \cite{mprox,cp}, and interior point algorithms \cite{nes_nem,wright_PD}. Which of these methods is practical depends on the specifics of the problem, such as the size and the cost of vector-matrix multiplications.

It is instructive to note that in the simplest setting of additive composite problems (Example~\ref{exa:add_comp}), the prox-linear method reduces to the popular proximal-gradient algorithm or ISTA \cite{smoothing_beckT}. For nonlinear least squares, the prox-linear method is a close variant of   Gauss-Newton. 

Recall that the step-size of the proximal point method provides a convenient stopping criteria, since it directly relates to the gradient of the Moreau envelope -- a smooth approximation of the objective function. Is there such an interpretation for the prox-linear method? This question is central, since termination criteria is not only used to stop the method but also to judge its efficiency and to compare against competing methods.

The answer is yes. Even though one can not evaluate the gradient $\|\nabla F_{\frac{1}{2\beta}}\|$ directly,
the scaled step-size of the prox-linear method $$\mathcal{G}(x):=\beta(x_{t+1}-x_t)$$ is a good surrogate \cite[Theorem 4.5]{prox_lin_paq}: 
$$\tfrac{1}{4} \|\nabla F_{\frac{1}{2\beta}}(x)\| \leq \|\mathcal{G}(x)\|\leq 3\|\nabla F_{\frac{1}{2\beta}}(x)\|.$$
In particular, the prox-linear method will find a point $x$ satisfying 
 $\|\nabla F_{\frac{1}{2\beta}}(x)\|^2\leq\varepsilon$ after at most $O\left(\frac{\beta(F(x_0)-\inf F)}{\varepsilon}\right)$ iterations. In the simplest setting when $g=0$ and $h(t)=t$, this rate reduces to the well-known convergence guarantee of gradient descent, which is black-box optimal for $C^1$-smooth nonconvex optimization \cite{grad_desc_opt}. 
 
It is worthwhile to note that a number of improvements to the basic prox-linear method were recently proposed. The authors of \cite{composite_cart} discuss trust region variants and their complexity guarantees, while \cite{duchi_ruan} propose stochastic extensions of the scheme and prove almost sure convergence. The paper \cite{prox_lin_paq} discusses overall complexity guarantees when the convex subproblems can only be solved by first-order methods, and proposes an inertial variant of the scheme whose convergence guarantees automatically adapt to the near-convexity of the problem.

\subsection{Local rapid convergence}

Under typical regularity conditions, the prox-linear method  exhibits the same types of rapid convergence guarantees as the proximal point method. I will illustrate with two intuitive and widely used regularity conditions, yielding local linear and quadratic convergence, respectively.

\begin{defn}[\cite{tilt}]{\rm
	A local minimizer $\bar x$ of $F$ is  {\em $\alpha$-tilt-stable}  if there exists $r>0$ such that the solution map
$$M: v\mapsto \argmin_{x\in B_r(\bar x)} \left\{ F(x)-\langle v,x \rangle\right\}$$
is $1/\alpha$-Lipschitz around $0$ with $M(0)=\bar x$. }
\end{defn}

This condition might seem unfamiliar to convex optimization specialist.
Though not obvious, tilt-stability is equivalent to a uniform quadratic growth property and  a subtle localization of strong convexity of $F$. See \cite{tilt_adrian} or \cite{Dima_Ng} for more details on these equivalences. Under the tilt-stability assumption, the prox-linear method initialized sufficiently close to $\bar x$ produces iterates that converge  at a linear rate $1-\alpha/\beta$.

The second regularity condition models sharp growth of the function around the minimizer. Let $S$ be the set of all stationary points of $F$, meaning $x$ lies in $S$ if and only if the directional derivative $F'(x;v)$ is nonnegative in every direction $v\in \R^d$.
\begin{defn}[\cite{weak_sharp}]{\rm
	A local minimizer $\bar x$ of $F$ is {\em sharp}  if there exists $\alpha>0$ and a neighborhood $\mathcal{X}$ of $\bar x$ such that 
$$F(x)\geq  F(\proj_S(x))+c\cdot \dist(x,S)\qquad\forall x\in \mathcal{X}.$$}
\end{defn}

Under the sharpness condition,  the prox-linear method initialized sufficiently close to $\bar x$ produces iterates that converge quadratically. 

For well-structured problems, one can hope to justify the two regularity conditions above under statistical assumptions. The recent work of Duchi and Ruan on the phase retrieval problem \cite{duchi_ruan_PR}  is an interesting recent example. Under mild statistical assumptions on the data generating mechanism,  sharpness  is assured with high probability. Therefore the  prox-linear method (and even subgradient methods \cite{proj_weak_dim}) converge rapidly, when initialized within a constant relative distance of an optimal solution.

\section{Catalyst acceleration}
The final example concerns inertial acceleration in convex optimization.  Setting the groundwork, consider a $\mu$-strongly convex function  $f$ with a  $\beta$-Lipschitz gradient map $x\mapsto \nabla f(x)$. Classically,  gradient descent  will find a point $x$ satisfying $f(x)-\min f<\varepsilon$ after at most $$O\left(\frac{\beta}{\mu}\ln(1/\varepsilon)\right)$$ iterations. Accelerated gradient methods, beginning with Nesterov~\cite{nest_orig}, equip the gradient descent method with an inertial correction. 
Such methods have the much lower complexity guarantee $$O\left(\sqrt{\frac{\beta}{\mu}}\ln(1/\varepsilon)\right),$$ which is optimal  within the first-order oracle model of computation \cite{complexity}. 
	
It is natural to ask which other methods, aside from gradient descent, can be ``accelerated''. For example, one may wish to accelerate coordinate descent or so-called variance reduced methods for finite sum problems; I will comment on the latter problem class shortly. 

One appealing strategy relies on the proximal point method. G\"{u}ler in \cite{gul_prox_acc} showed that the proximal point method itself can be equipped with inertial steps leading to improved convergence guarantees. Building on this work, Lin, Mairal, and Harchaoui \cite{catalyst} explained how to derive the {\em total} complexity guarantees for an inexact accelerated proximal point method that take into account the cost of applying
 an arbitrary linearly convergent algorithm $\mathcal{M}$ to the subproblems. Their {\em Catalyst acceleration} framework is summarized in Algorithm~\ref{alg:catalyst}.

\begin{algorithm}
	\KwData{$x_0\in \R^d$, $\kappa>0$, algorithm $\mathcal{M}$}
Set $q= \mu/(\mu+\kappa)$, $\alpha_0=\sqrt{q}$, and $y_0=x_0$\;
	\For{t=0,\ldots,T}{
	Use $\mathcal{M}$ to approximately solve:
\begin{equation}\label{eqn:prox_subprob}
	x_t\approx\argmin_{x\in \R^d} \left\{F(x)+\frac{\kappa}{2}\|x-y_{t-1}\|^2\right\}.\;
\end{equation}

Compute $\alpha_t\in (0,1)$ from the equation $$\alpha_t^2=(1-\alpha_t)\alpha_{t-1}^2+q\alpha_t.\;$$

Compute:
\begin{align*}
\beta_t&=\frac{\alpha_{t-1}(1-\alpha_{t-1})}{\alpha_{t-1}^2+\alpha_t},\\
y_t&=x_t+\beta_t(x_t-x_{t-1}). 
\end{align*}
}
\caption{Catalyst Acceleration}\label{alg:catalyst}
\end{algorithm}

To state the guarantees of this method, suppose that $\mathcal{M}$ converges on the proximal subproblem in function value at a linear rate $1-\tau\in (0,1)$. Then a simple termination policy on the subproblems \eqref{eqn:prox_subprob} yields an algorithm with overall complexity 
\begin{equation}\label{eqn:compl}
\widetilde{O}\left(\frac{\sqrt{\mu+\kappa}}{\tau \sqrt{\mu}}\ln(1/\varepsilon)\right).
\end{equation}
That is, the expression \eqref{eqn:compl} describes the maximal number of iterations of $\mathcal{M}$ used by Algorithm~\ref{alg:catalyst} until it finds a point $x$ satisfying $f(x)-\inf f\leq \varepsilon$.
Typically $\tau$ depends on $\kappa$; therefore the best choice of $\kappa$ is the one that minimizes the ratio $\frac{\sqrt{\mu+\kappa}}{\tau \sqrt{\mu}}$.

The main motivation for the Catalyst framework, and its most potent application, is the regularized Empirical Risk Minimization (ERM) problem:
$$\min_{x\in \R^d} f(x):=\frac{1}{m}\sum_{i=1}^m f_i(x)+g(x).$$
Such large-finite sum problems are ubiquitous in machine learning and high-dimensional statistics, where each function $f_i$ typically models a misfit between predicted and observed data while $g$ promotes some low dimensional structure on $x$, such as sparsity or low-rank.

Assume that $f$ is $\mu$-strongly convex and each individual $f_i$ is $ C^1$-smooth with $\beta$-Lipschitz gradient. Since $m$ is assumed to be huge, the complexity of numerical methods is best measured in terms of the total number of individual gradient evaluations $\nabla f_i$. In particular, fast gradient methods have the worst-case complexity $$O\left(m\sqrt{\frac{\beta}{\mu}}\ln(1/\varepsilon)\right),$$ since each iteration requires evaluation of all the individual gradients $\{\nabla f_i(x)\}_{i=1}^m$. Variance reduced algorithms, such as SAG \cite{sag}, SAGA \cite{SAGA2}, SDCA \cite{sdca}, SMART \cite{smart_davis}, SVRG \cite{svrg,prox_SVRG}, FINITO \cite{finito}, and MISO \cite{miso,catalyst}, aim to improve the dependence on $m$. In their raw form, all of these methods exhibit a similar complexity $$O\left(\left(m+\frac{\beta}{\mu}\right)\ln(1/\varepsilon)\right),$$ in expectation, and differ only in storage requirements and in whether one needs to know explicitly the strong convexity constant. 

It was a long standing open question to determine if the dependence on $\beta/\mu$ can be improved. This is not quite possible in full generality, and instead one should expect a rate of the form
 $$O\left(\left(m+\sqrt{m\frac{\beta}{\mu}}\right)\ln(1/\varepsilon)\right).$$ 
Indeed, such a rate would be optimal in an appropriate oracle model of complexity \cite{NIPS2016_6058,yosi,AgB,conjugategradient}. Thus acceleration for ERM problems is only beneficial in the setting $m< \beta/\mu$.

  Early examples for specific algorithms are the accelerated SDCA \cite{accsdca} and RPDG \cite{conjugategradient}.\footnote{Here, I am ignoring logarithmic terms in the convergence rate.} The accelerated SDCA, in particular, uses a specialized proximal-point construction and  was the motivation for the Catalyst framework. Catalyst generic acceleration allows to accelerate all of the variance reduced methods above in a single conceptually transparent framework.  It is worth noting that the first direct accelerated variance reduced methods for ERM problems were recently proposed in \cite{accsvrg,NIPS2016_6154}. 
  
  In contrast to the convex setting, the role of inertia for nonconvex problems is not nearly as well understood. In particular, gradient descent is black-box optimal for $C^1$-smooth nonconvex minimization \cite{grad_desc_opt}, and therefore inertia can not help in the worst case. On the other hand, the recent paper \cite{pmlr-v70-carmon17a} presents a first-order method for minimizing $C^2$ and $C^3$ smooth functions that is provably faster than gradient descent. At its core, their algorithm also combines inertia with the proximal point method.
  For a partial extension of the Catalyst framework to weakly convex problems, see \cite{catalyst_2}. 

\section{Conclusion}
The proximal point method has long been ingrained in the foundations of optimization.  Recent  progress in large scale computing has shown that the proximal point method is not only conceptual, but can guide methodology. Though direct methods are usually preferable, proximally guided algorithms can be equally effective and often lead to more easily interpretable numerical methods.
In this article, I outlined three examples of this viewpoint, where the proximal-point method guides both the design and analysis of numerical methods.  
	\bigskip
	
\noindent{\bf Acknowledgments.}  The author thanks Damek Davis, John Duchi, and Zaid Harchaoui for their helpful comments on an early draft of the article. Research of Drusvyatskiy is supported by the AFOSR YIP award FA9550-15-1-0237 and by the NSF DMS   1651851 and CCF 1740551 awards.




\bibliographystyle{plain}
\bibliography{bibliography}

\def\cfac#1{\ifmmode\setbox7\hbox{$\accent"5E#1$}\else
  \setbox7\hbox{\accent"5E#1}\penalty 10000\relax\fi\raise 1\ht7
  \hbox{\lower1.15ex\hbox to 1\wd7{\hss\accent"13\hss}}\penalty 10000
  \hskip-1\wd7\penalty 10000\box7}
\begin{thebibliography}{10}

\bibitem{abbe_band}
E.~Abbe, A.S. Bandeira, A.~Bracher, and A.~Singer.
\newblock Decoding binary node labels from censored edge measurements: phase
  transition and efficient recovery.
\newblock {\em IEEE Trans. Network Sci. Eng.}, 1(1):10--22, 2014.

\bibitem{AgB}
A.~Agarwal and L.~Bottou.
\newblock A lower bound for the optimization of finite sums.
\newblock In {\em Proceedings of the 32nd International Conference on Machine
  Learning, {ICML} 2015, Lille, France, 6-11 July 2015}, pages 78--86, 2015.

\bibitem{accsvrg}
Z.~Allen-Zhu.
\newblock Katyusha: The first direct acceleration of stochastic gradient
  methods.
\newblock {\em Preprint arXiv:1603.05953 (version 5)}, 2016.

\bibitem{yosi}
Y.~Arjevani.
\newblock Limitations on variance-reduction and acceleration schemes for finite
  sums optimization.
\newblock In I.~Guyon, U.~V. Luxburg, S.~Bengio, H.~Wallach, R.~Fergus,
  S.~Vishwanathan, and R.~Garnett, editors, {\em Advances in Neural Information
  Processing Systems 30}, pages 3543--3552. Curran Associates, Inc., 2017.

\bibitem{ban_boum}
A.S. Bandeira, N.~Boumal, and V.~Voroninski.
\newblock On the low-rank approach for semidefinite programs arising in
  synchronization and community detection.
\newblock In {\em Proceedings of the 29th Conference on Learning Theory, {COLT}
  2016, New York, USA, June 23-26, 2016}, pages 361--382, 2016.

\bibitem{jordan}
P.L. Bartlett, M.I. Jordan, and J.D. McAuliffe.
\newblock Convexity, classification, and risk bounds.
\newblock {\em J. Amer. Statist. Assoc.}, 101(473):138--156, 2006.

\bibitem{BC_book}
H.H. Bauschke and P.L. Combettes.
\newblock {\em Convex analysis and monotone operator theory in {H}ilbert
  spaces}.
\newblock CMS Books in Mathematics/Ouvrages de Math\'ematiques de la SMC.
  Springer, Cham, second edition, 2017.
\newblock With a foreword by H\'edy Attouch.

\bibitem{smoothing_beckT}
A.~Beck and M.~Teboulle.
\newblock Smoothing and first order methods: a unified framework.
\newblock {\em SIAM J. Optim.}, 22(2):557--580, 2012.

\bibitem{BB}
L.~Bottou and O.~Bousquet.
\newblock The tradeoffs of large scale learning.
\newblock In {\em Advances in Neural Information Processing Systems}, pages
  161--168, 2008.

\bibitem{weak_sharp}
J.V. Burke and M.C. Ferris.
\newblock Weak sharp minima in mathematical programming.
\newblock {\em SIAM J. Control Optim.}, 31(5):1340--1359, 1993.

\bibitem{quad_conv}
J.V. Burke and M.C. Ferris.
\newblock A {G}auss-{N}ewton method for convex composite optimization.
\newblock {\em Math. Programming}, 71(2, Ser. A):179--194, 1995.

\bibitem{rob_cand}
E.J. Cand{\`e}s, X.~Li, Y.~Ma, and J.~Wright.
\newblock Robust principal component analysis?
\newblock {\em J. ACM}, 58(3):Art. 11, 37, 2011.

\bibitem{wirt_flow}
E.J. Cand\`es, X.~Li, and M.~Soltanolkotabi.
\newblock Phase retrieval via {W}irtinger flow: theory and algorithms.
\newblock {\em IEEE Trans. Inform. Theory}, 61(4):1985--2007, 2015.

\bibitem{pmlr-v70-carmon17a}
Y.~Carmon, J.C. Duchi, O.~Hinder, and A.~Sidford.
\newblock ``{C}onvex until proven guilty'': Dimension-free acceleration of
  gradient descent on non-convex functions.
\newblock In {\em Proceedings of the 34th International Conference on Machine
  Learning}, volume~70, pages 654--663, 2017.

\bibitem{grad_desc_opt}
Y.~Carmon, J.C. Duchi, O.~Hinder, and A.~Sidford.
\newblock Lower bounds for finding stationary points i.
\newblock {\em Preprint arXiv:1710.11606}, 2017.

\bibitem{composite_cart}
C.~Cartis, N.I.M. Gould, and P.L. Toint.
\newblock On the evaluation complexity of composite function minimization with
  applications to nonconvex nonlinear programming.
\newblock {\em SIAM J. Optim.}, 21(4):1721--1739, 2011.

\bibitem{cp}
A.~Chambolle and T.~Pock.
\newblock A first-order primal-dual algorithm for convex problems with
  applications to imaging.
\newblock {\em J. Math. Imaging Vision}, 40(1):120--145, 2011.

\bibitem{chand}
V.~Chandrasekaran, S.~Sanghavi, P.~A. Parrilo, and A.S. Willsky.
\newblock Rank-sparsity incoherence for matrix decomposition.
\newblock {\em SIAM J. Optim.}, 21(2):572--596, 2011.

\bibitem{rand_quad}
Y.~Chen and E.J. Cand\`es.
\newblock Solving random quadratic systems of equations is nearly as easy as
  solving linear systems.
\newblock {\em Comm. Pure Appl. Math.}, 70(5):822--883, 2017.

\bibitem{smart_davis}
D.~Davis.
\newblock {SMART}: The stochastic monotone aggregated root-finding algorithm.
\newblock {\em Preprint arXiv:1601.00698}, 2016.

\bibitem{proj_weak_dim}
D.~Davis, D.~Drusvyatskiy, and C.~Paquette.
\newblock The nonsmooth landscape of phase retrieval.
\newblock {\em Preprint arXiv:1711.03247}, 2017.

\bibitem{prox_guide_subgrad}
D.~Davis and B.~Grimmer.
\newblock Proximally guided stochastic sbgradient method for nonsmooth,
  nonconvex problems.
\newblock {\em Preprint, arXiv:1707.03505}, 2017.

\bibitem{NIPS2016_6154}
A.~Defazio.
\newblock A simple practical accelerated method for finite sums.
\newblock In D.~D. Lee, M.~Sugiyama, U.~V. Luxburg, I.~Guyon, and R.~Garnett,
  editors, {\em Advances in Neural Information Processing Systems 29}, pages
  676--684. Curran Associates, Inc., 2016.

\bibitem{SAGA2}
A.~Defazio, F.~Bach, and S.~Lacoste-Julien.
\newblock {SAGA}: A fast incremental gradient method with support for
  non-strongly convex composite objectives.
\newblock In Z.~Ghahramani, M.~Welling, C.~Cortes, N.~D. Lawrence, and K.~Q.
  Weinberger, editors, {\em Advances in Neural Information Processing Systems
  27}, pages 1646--1654. Curran Associates, Inc., 2014.

\bibitem{finito}
A.~Defazio, J.~Domke, and T.S. Caetano.
\newblock Finito: A faster, permutable incremental gradient method for big data
  problems.
\newblock In {\em ICML}, pages 1125--1133, 2014.

\bibitem{tilt_adrian}
D.~Drusvyatskiy and A.S. Lewis.
\newblock Tilt stability, uniform quadratic growth, and strong metric
  regularity of the subdifferential.
\newblock {\em SIAM J. Optim.}, 23(1):256--267, 2013.

\bibitem{prox_error}
D.~Drusvyatskiy and A.S. Lewis.
\newblock Error bounds, quadratic growth, and linear convergence of proximal
  methods.
\newblock {\em To appear in Math. Oper. Res., arXiv:1602.06661}, 2016.

\bibitem{Dima_Ng}
D.~Drusvyatskiy, B.S. Mordukhovich, and T.T.A. Nghia.
\newblock Second-order growth, tilt-stability, and metric regularity of the
  subdifferential.
\newblock {\em J. Convex Anal.}, 21(4).

\bibitem{prox_lin_paq}
D.~Drusvyatskiy and C.~Paquette.
\newblock Efficiency of minimizing compositions of convex functions and smooth
  maps.
\newblock {\em Preprint, arXiv:1605.00125}, 2016.

\bibitem{duchi_ruan_PR}
J.C. Duchi and F.~Ruan.
\newblock Solving (most) of a set of quadratic equalities: Composite
  optimization for robust phase retrieval.
\newblock {\em Preprint arXiv:1705.02356}, 2017.

\bibitem{duchi_ruan}
J.C. Duchi and F.~Ruan.
\newblock Stochastic methods for composite optimization problems.
\newblock {\em Preprint arXiv:1703.08570}, 2017.

\bibitem{eM}
Y.C. Eldar and S.~Mendelson.
\newblock Phase retrieval: stability and recovery guarantees.
\newblock {\em Appl. Comput. Harmon. Anal.}, 36(3):473--494, 2014.

\bibitem{gl_stoch}
S.~Ghadimi and G.~Lan.
\newblock Stochastic first- and zeroth-order methods for nonconvex stochastic
  programming.
\newblock {\em SIAM J. Optim.}, 23(4):2341--2368, 2013.

\bibitem{gul_prox_acc}
O.~G\"uler.
\newblock New proximal point algorithms for convex minimization.
\newblock {\em SIAM J. Optim.}, 2(4):649--664, 1992.

\bibitem{hazan_subgrad}
E.~Hazan and S.~Kale.
\newblock Beyond the regret minimization barrier: an optimal algorithm for
  stochastic strongly-convex optimization.
\newblock In Sham~M. Kakade and Ulrike von Luxburg, editors, {\em Proceedings
  of the 24th Annual Conference on Learning Theory}, volume~19 of {\em
  Proceedings of Machine Learning Research}, pages 421--436, Budapest, Hungary,
  09--11 Jun 2011. PMLR.

\bibitem{svrg}
R.~Johnson and T.~Zhang.
\newblock Accelerating stochastic gradient descent using predictive variance
  reduction.
\newblock In {\em Proceedings of the 26th International Conference on Neural
  Information Processing Systems}, NIPS'13, pages 315--323, USA, 2013. Curran
  Associates Inc.

\bibitem{MR3353214}
A.~Juditsky and Y.~Nesterov.
\newblock Deterministic and stochastic primal-dual subgradient algorithms for
  uniformly convex minimization.
\newblock {\em Stoch. Syst.}, 4(1):44--80, 2014.

\bibitem{a_simp_app}
S.~Lacoste-Julien, M.~Schmidt, and F.~Bach.
\newblock A simpler approach to obtaining an ${O}(1/t)$ convergence rate for
  the projected stochastic subgradient method.
\newblock {\em arxiv arXiv:1212.2002}, 2012.

\bibitem{conjugategradient}
G.~Lan.
\newblock An optimal randomized incremental gradient method.
\newblock {\em arXiv:1507.02000}, 2015.

\bibitem{bundle}
C.~Lemarechal, J.-J. Strodiot, and A.~Bihain.
\newblock On a bundle algorithm for nonsmooth optimization.
\newblock In {\em Nonlinear programming, 4 ({M}adison, {W}is., 1980)}, pages
  245--282. Academic Press, New York-London, 1981.

\bibitem{prox}
A.S. Lewis and S.J. Wright.
\newblock A proximal method for composite minimization.
\newblock {\em Math. Program.}, pages 1--46, 2015.

\bibitem{catalyst}
H.~Lin, J.~Mairal, and Z.~Harchaoui.
\newblock A universal catalyst for first-order optimization.
\newblock In {\em Advances in Neural Information Processing Systems}, pages
  3366--3374, 2015.

\bibitem{luke_news_views}
R.~Luke.
\newblock {Phase Retrieval, What's New?}
\newblock {\em SIAG/OPT Views and News}, 25(1), 2017.

\bibitem{miso}
J.~Mairal.
\newblock Incremental majorization-minimization optimization with application
  to large-scale machine learning.
\newblock {\em SIAM Journal on Optimization}, 25(2):829--855, 2015.

\bibitem{mprox}
A.~Nemirovski.
\newblock Prox-method with rate of convergence {$O(1/t)$} for variational
  inequalities with {L}ipschitz continuous monotone operators and smooth
  convex-concave saddle point problems.
\newblock {\em SIAM J. Optim.}, 15(1):229--251, 2004.

\bibitem{latest_subgrad}
A.~Nemirovski, A.~Juditsky, G.~Lan, and A.~Shapiro.
\newblock Robust stochastic approximation approach to stochastic programming.
\newblock {\em SIAM J. Optim.}, 19(4):1574--1609, 2008.

\bibitem{complexity}
A.S. Nemirovsky and D.B. Yudin.
\newblock {\em Problem complexity and method efficiency in optimization}.
\newblock A Wiley-Interscience Publication. John Wiley \& Sons, Inc., New York,
  1983.
\newblock Translated from the Russian and with a preface by E. R. Dawson,
  Wiley-Interscience Series in Discrete Mathematics.

\bibitem{nes_nem}
Y.~Nesterov and A.~Nemirovskii.
\newblock {\em Interior-point polynomial algorithms in convex programming},
  volume~13 of {\em SIAM Studies in Applied Mathematics}.
\newblock Society for Industrial and Applied Mathematics (SIAM), Philadelphia,
  PA, 1994.

\bibitem{nest_orig}
Yu. Nesterov.
\newblock A method for solving the convex programming problem with convergence
  rate {$O(1/k^{2})$}.
\newblock {\em Dokl. Akad. Nauk SSSR}, 269(3):543--547, 1983.

\bibitem{smooth_min_nonsmooth}
Yu. Nesterov.
\newblock Smooth minimization of non-smooth functions.
\newblock {\em Math. Program.}, 103(1, Ser. A):127--152, 2005.

\bibitem{nest_GN}
Yu. Nesterov.
\newblock Modified {G}auss-{N}ewton scheme with worst case guarantees for
  global performance.
\newblock {\em Optim. Methods Softw.}, 22(3):469--483, 2007.

\bibitem{nest_conv_comp}
Yu. Nesterov.
\newblock Gradient methods for minimizing composite functions.
\newblock {\em Math. Program.}, 140(1, Ser. B):125--161, 2013.

\bibitem{catalyst_2}
C.~Paquette, H.~Lin, D.~Drusvyatskiy, J.~Mairal, and Z.~Harchaoui.
\newblock Catalyst acceleration for gradient-based non-convex optimization.
\newblock {\em Preprint arXiv:1703.10993}, 2017.

\bibitem{tilt}
R.A. Poliquin and R.T. Rockafellar.
\newblock Tilt stability of a local minimum.
\newblock {\em SIAM J. Optim.}, 8(2):287--299 (electronic), 1998.

\bibitem{stochave}
B.T. Polyak and A.B. Juditsky.
\newblock Acceleration of stochastic approximation by averaging.
\newblock {\em SIAM J. Control Optim.}, 30(4):838--855, 1992.

\bibitem{Rakhlin_subgrad}
A.~Rakhlin, O.~Shamir, and K.~Sridharan.
\newblock Making gradient descent optimal for strongly convex stochastic
  optimization.
\newblock In {\em Proceedings of the 29th International Coference on
  International Conference on Machine Learning}, ICML'12, pages 1571--1578,
  USA, 2012. Omnipress.

\bibitem{rob_mon}
H.~Robbins and S.~Monro.
\newblock A stochastic approximation method.
\newblock {\em Ann. Math. Statistics}, 22:400--407, 1951.

\bibitem{mon_rock}
R.T. Rockafellar.
\newblock Monotone operators and the proximal point algorithm.
\newblock {\em SIAM J. Control Optimization}, 14(5):877--898, 1976.

\bibitem{sag}
M.~Schmidt, N.~Le Roux, and F.~Bach.
\newblock Minimizing finite sums with the stochastic average gradient.
\newblock {\em arXiv:1309.2388}, 2013.

\bibitem{sdca}
S.~Shalev-Shwartz and T.~Zhang.
\newblock Proximal stochastic dual coordinate ascent.
\newblock {\em arXiv:1211.2717}, 2012.

\bibitem{accsdca}
S.~Shalev-Shwartz and T.~Zhang.
\newblock Accelerated proximal stochastic dual coordinate ascent for
  regularized loss minimization.
\newblock {\em Mathematical Programming}, 2015.

\bibitem{ang_sing}
A.~Singer.
\newblock Angular synchronization by eigenvectors and semidefinite programming.
\newblock {\em Appl. Comput. Harmon. Anal.}, 30(1):20--36, 2011.

\bibitem{phase_nonconv}
J.~Sun, Q.~Qu, and J.~Wright.
\newblock A geometric analysis of phase retrieval.
\newblock {\em To appear in Found. Comp. Math., arXiv:1602.06664}, 2017.

\bibitem{NIPS2016_6058}
B.E. Woodworth and N.~Srebro.
\newblock Tight complexity bounds for optimizing composite objectives.
\newblock In D.~D. Lee, M.~Sugiyama, U.~V. Luxburg, I.~Guyon, and R.~Garnett,
  editors, {\em Advances in Neural Information Processing Systems 29}, pages
  3639--3647. Curran Associates, Inc., 2016.

\bibitem{wright_PD}
S.J. Wright.
\newblock {\em Primal-dual interior-point methods}.
\newblock Society for Industrial and Applied Mathematics (SIAM), Philadelphia,
  PA, 1997.

\bibitem{prox_SVRG}
L.~Xiao and T.~Zhang.
\newblock A proximal stochastic gradient method with progressive variance
  reduction.
\newblock {\em SIAM J. Optim.}, 24(4):2057--2075, 2014.

\end{thebibliography}

\end{document}